\documentclass[journal, 3p]{elsarticle}

\usepackage[utf8]{inputenc}
\usepackage[T1]{fontenc}
\usepackage{lmodern}

\makeatletter
\def\ps@pprintTitle{%
	\let\@oddhead\@empty
	\let\@evenhead\@empty
	\let\@oddfoot\@empty
	\let\@evenfoot\@oddfoot
}
\makeatother

\usepackage{algorithm}
\usepackage{algpseudocode}
\usepackage{comment}
\newcommand{\dA}{\Delta} 
\newcommand{\A}{A}   
\newcommand{\dAi}{\Delta_i}

\newcommand{\sign}{\mbox{sign}}

\newcommand{\R}{\mathbb{R}}
\newcommand{\U}{\mathcal{U}}
\newcommand{\ith}{ ^{\text{th}} }
\newcommand{\constraint}{ {\|\dA\|_{\infty} \le \delta} }

\newcommand{\oner}{\mathbf{1}}

\usepackage{url}

\usepackage{graphicx, caption}
\usepackage{amsmath,amsfonts,amssymb, amsthm}
\usepackage[usenames,dvipsnames,svgnames]{xcolor} 
\usepackage{booktabs}
\usepackage{siunitx} 
\usepackage{etoolbox}

\newtheorem{theorem}{Theorem}
\newtheorem{corollary}{Corollary}

\begin{document}

\begin{frontmatter}
 \title{Robust Least Squares for Quantized Data Matrices}

\author{Stephen Becker}
\ead{stephen.becker@colorado.edu}

\cortext[cor1]{Corresponding author}
\author{Richard Clancy\corref{cor1}}
\ead{richard.clancy@colorado.edu}

\address{Department of Applied Mathematics at the University of Colorado, Boulder, CO 80309}


\begin{abstract}
	In this paper we formulate and solve a robust least squares problem for a system of linear equations subject to quantization error in the data matrix. Ordinary least squares fails to consider uncertainty in the operator, modeling all noise in the observed signal. Total least squares accounts for uncertainty in the data matrix, but necessarily increases the condition number of the operator compared to ordinary least squares. Tikhonov regularization or ridge regression is frequently employed to combat ill-conditioning, but requires parameter tuning which presents a host of challenges and places strong assumptions on parameter prior distributions. The proposed method also requires selection of a parameter, but it can be chosen in a natural way, e.g., a matrix rounded to the 4th digit uses an uncertainty bounding parameter of \num{0.5e-4}. We show here that our robust method is theoretically appropriate, tractable, and performs favorably against ordinary and total least squares. 
\end{abstract}

\begin{keyword}
	\noindent Least squares approximation \sep Quantization error \sep Ridge regression \sep Robust optimization \sep Subgradient methods \sep Tikhonov regularization \sep Total least squares
\end{keyword}

\end{frontmatter}

\section{Introduction}
The primary goal of this paper is to recover unknown parameters from a noisy observation and an uncertain linear operator. In particular, our mathematical model for robust least squares is
\begin{equation}\label{eq:main}
	\min_x \, \left\{ \max_{ \dA \in \,  \U }\; \| (\A+\dA)x -b \|^2 \right\} 
\end{equation}
which we refer to as our robust optimization (RO) problem. The Euclidean norm is denoted by $\| \cdot \|$ and $\mathcal{U}$ is the uncertainty set from which perturbations in $A$ are drawn.  

The above RO formulation is motivated by two situations. In both cases, we let $\bar A$ and $\bar x$ represent the \textit{true and unknown} data matrix and parameter vector, respectively. We model $b = \bar A \bar x + \eta$, with $\eta$ i.i.d.\ Gaussian, and we only have knowledge of $A = \bar A - \dA$. We don't know $\dA$ explicitly but can make inferences based on the problem. In the first situation, we consider a data matrix subject to quantization or round-off error. Suppose the observed matrix $A$ has elements rounded to the hundredth place. Our uncertainty set can be written as $\U= \{ \dA \in \R^{m \times n } : \| \dA\|_{\infty} \le \delta \} $ with $\delta = 0.005$. If $A_{i,j} = 0.540$, then we know the true $\bar A_{i,j} \in (0.535,\, 0.545]$, hence $\dA_{i,j} \in (-0.005, 0.005]$. The norm $\| \cdot \|_{\infty}$ takes the maximum absolute value of any element in the matrix.

The second problem considers a data matrix with uncertainty proportional to the magnitude of each entry, i.e. $\U= \{ \dA \in \R^{m \times n} : \dA_{i,j} \in (-p |\A_{i,j}|, p |\A_{i,j}|  ]   \}$. Here, $p$ denotes a proportionality constant. Data subject to $\pm 1 \%$ uncertainty would have $p = 0.01$. 
The two cases cover the effects of finite-precision in fixed and floating point representations, respectively. In both problems, the uncertainty sets are specified element-wise allowing us to decouple along rows. 

Due to limitations in both ordinary (OLS) and total least squares (TLS), the signal processing community has sought alternatives when operator uncertainty is present. In \cite{wiesel2008linear} and \cite{zhu2014maximum}, the authors derive maximum likelihood estimators (MLE) and Cram\'er Rao bounds for problems with data matrices subject to Gaussian noise under two different models; errors-in-variables and a random variable model. The setup for \cite{zhu2014maximum} only has access to sign measurements which is more restrictive than we consider here. Both papers focus primarily on the case where variance in $A$ and $b$ are known; 
in contrast, we make no distributional assumption on $A$ (it could be stochastic, deterministic, or adversarial). In the absence of knowing the variance, \cite{wiesel2008linear} shows their estimator reduces to the OLS solution. The problem is treated through an approximate message passing framework in \cite{zhu2020vector} for more general, structured perturbations. Our use case fits under this umbrella but they impose a sparsity inducing prior. Although our method performs well with sparse solutions, we do not assume it.

Note that the model in Eq.~\ref{eq:main} differs from other robust least squares problems considered in \cite{CaramanisLasso}, \cite{SmolaRobust07} and classic papers from the late 1990's \cite{el1997robust} and early 2000's \cite{goldfarb2003robust}. Those works usually made special assumptions like the constraint and objective norms match (e.g. minimize $\|x\|_p$ subject to $\|Ax- b\|_p \le v$), column-wise separability, or ellipsoidal uncertainty sets. The work in \cite{SmolaRobust07} is similar regarding row-wise separability, but they impose a hard constraint on the 2-norm of the solution, i.e., $\|x\| \le k$ for $k \in \R$. For an overview of robust optimization, see \cite{BertsimasCaramanisRO} and \cite{gorissen2015practical} with the latter focusing on implementation.

In this paper, we consider box constraints over entries of the data matrix. Our main contribution is the formulation of a robust objective to handle quantization error in least squares problems and the presentation of methods to solve it. Although the proposed method requires selection of an uncertainty parameter, it is chosen in a natural and theoretically appropriate way based on the observed extent of quantization. This is in contrast to ridge regression and MLE based methods that require involved parameter tuning or \textit{a priori} knowledge of the probability distributions from which model uncertainty is drawn. We anticipate our method to be most effective under moderate to heavy quantization where fidelity loss is greater than $0.1\%$.

\subsection{Limitations of Ordinary and Total Least Squares} \label{sec:ols}
The ordinary least squares (OLS) problem seeks a vector $x$ to minimize the residual given by
\begin{equation} \label{eq:OLS}
    \min_x \|\A x-b\|^2.
\end{equation}
We focus our attention on the over-determined case where $m>n$ and further assume that $A$ is full rank. It is well known that the closed form solution to (\ref{eq:OLS}) is given by 
\begin{equation} \label{eq:OLS_solution}
    \hat x_\text{OLS} = (\A^T\A)^{-1}A^T b.
\end{equation}
A key assumption for OLS is that $A$ is known with certainty and $b$ is subject to additive noise. 
The OLS solution \eqref{eq:OLS_solution} is the MLE for the model $Ax-b = \eta \sim \mathcal{N}(0, \sigma^2 I)$. 

In practice, it is uncommon to know $A$ precisely. Typical causes of uncertainty are sampling error, measurement error, human error, modeling error, or rounding error. There were attempts at addressing this model limitation in \cite{hodge1972data} but these relied on small magnitude errors to use a truncated series approximation for an inverse. Total least square (TLS) was developed in response to this lack of symmetry in uncertainty. Rather than isolating noise to the observed signal or right hand side, the TLS model allows for an uncertain operator and is given by $(A + \dA) \bar x = b + \eta$ where $A$ and $b$ are observed with $\dA$ and $\eta$ a random matrix and vector, respectively. We assume that $\dA$ and $\eta$ have been scaled to have the same variance; see \cite{vanhuffel2007tls} for a modern overview of the topic. The TLS solutions, $\hat x_{\text{TLS}}$, solves 
\begin{equation} \label{eq:TLS}
    \min_{x,\dA,\eta}  \quad \big\| [\dA, \  \ \eta] \big\|_F \quad \text{subject to} \quad  (\A + \dA)x = b + \eta 
\end{equation}
where $[\dA, \ \eta] \in \R^{m \times (n+1)}$, $A$ and $b$ are the observed matrix/signal pair, and $\| \cdot \|_F$ is the Frobenius norm. It was shown in \cite{golub1980tls} that (\ref{eq:TLS}) can be solved via a singular value decomposition (SVD) with the closed form solution of
\begin{equation} \label{eq:tls_solution}
    \hat x_\text{TLS} = (\A^T \A - \sigma_{n+1}^2 I)^{-1}\A^T b
\end{equation}
where $\sigma_{n+1} \in \R$ being the smallest singular value of the augmented matrix $[A, \ b] \in \R^{m \times (n+1)}$. Similar to OLS, the TLS solution yields a MLE for the model of Gaussian noise in $A$ and $b$. It should be noted that $(\A^T \A - \sigma_{n+1}^2 I)$ is necessarily worse conditioned than $\A^T \A$. Since $\A^T \A$ is positive definite by virtue of $\A$ being full rank, all eigenvalues of $\A^T \A$ are shifted closer to zero by the amount of $\sigma_{n+1}^2$. For a small spectral gap between $\sigma_n$ and $\sigma_{n+1}$, the matrix will be close to singular making solutions extremely sensitive to perturbations in $A$. 

\begin{figure}  
	\centerline{\includegraphics[width=.4 \textwidth]{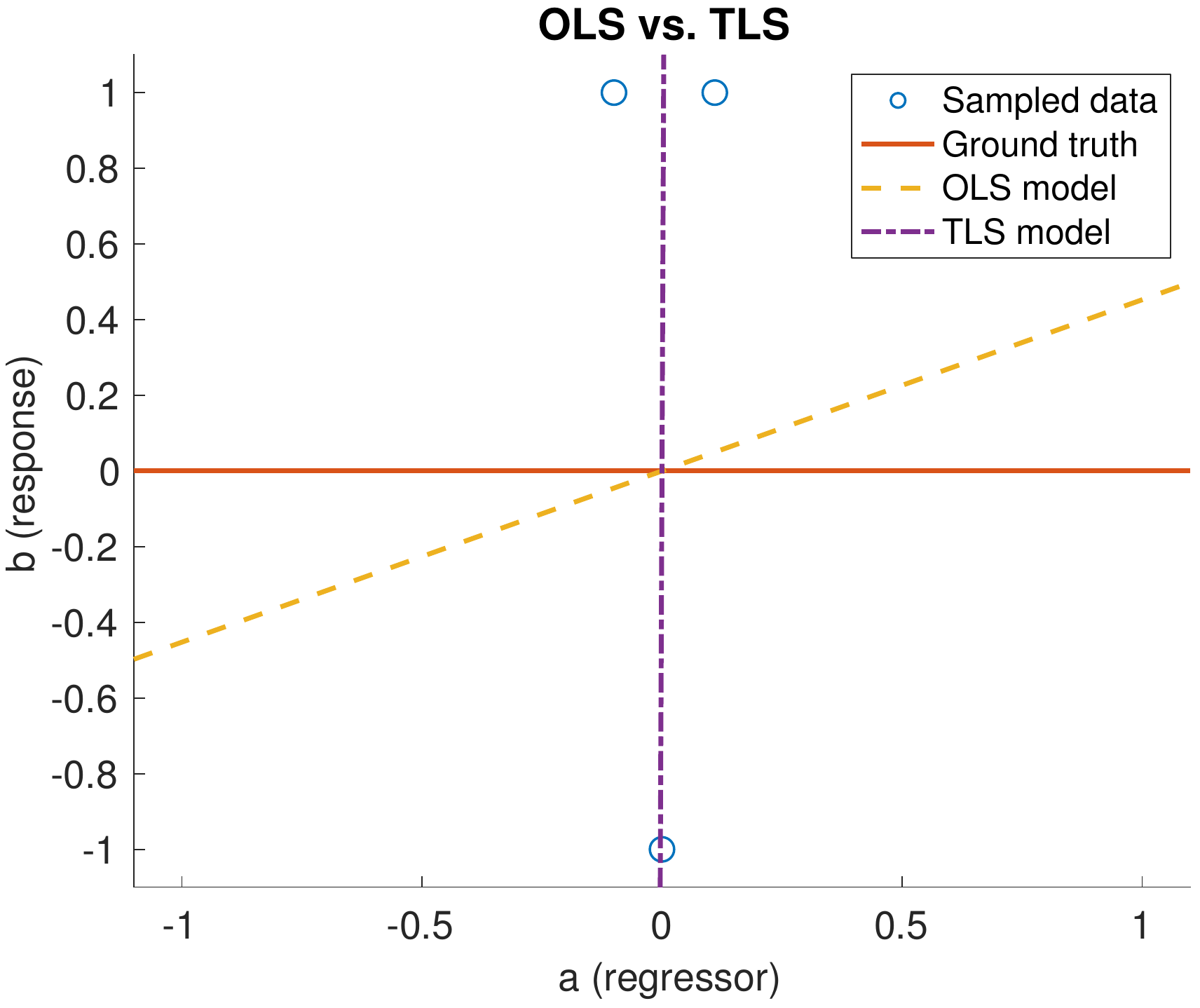}} \
	\caption{Instance of data sampled from a constant function, $\bar b = 0$ subject to additive noise. That is, $b = \bar b + \delta = \delta$. Solid line indicates ground truth model ($\bar b = 0$). Dashed line shows OLS model ($b = a \cdot \hat x_{\text{OLS}}$). TLS model is dot-dashed line ($b = a \cdot \hat x_{\text{TLS}}$). Since TLS minimizes orthogonal distance, the model over-fits data.}
	\label{fig:ols-vs-tls}
\end{figure}

This can be understood intuitively; uncertainty in $A$ permits additional degrees of freedom allowing the model to ``fit'' noise. To illustrate, consider the simple linear regression problem. We have access to several regressor / response pairs $(a, b) \in \R^2$. Suppose our data is generated from the model $b = a \cdot 0 + \delta$ with $\delta \sim Uniform(\{-1,1\})$ (the zero function plus discrete and uniform noise). We'd like to recover the true slope parameter, $x= 0$. In this instance, our three samples are given by $(-0.10,\ 1.00)$, $(0.00,\ -1.00)$, and $(0.11,\ 1.00)$. The OLS solution returns a slope of $\hat x_{\text{OLS}}\approx 0.45$ whereas TLS gives $\hat x_{\text{TLS}} \approx 297.79$. This is shown graphically in Figure \ref{fig:ols-vs-tls}. TLS's ability to vary the operator results in extreme sensitivity to data provided. Although our example appears exceptional, this behavior is typical of TLS. Because of this sensitivity, a number of alternatives have been studied for uncertain operators as mentioned above. A more traditional approach for addressing the ill-conditioning encountered in TLS is ridge regression or Tikhonov regularization which we consider below.

\subsection{The Trouble with Tikhonov Regularization} \label{sec:Tikhonov}
Poor conditioning as observed in TLS is often combated with Tikhonov regularization or ridge regression (RR) whereby solutions with large norms are penalized via an $\ell_2$ regularization term. The solution to the ridge regression problem solves
\begin{equation} \label{eq:RR}
    \min_x \bigg\{\|\A x-b\|^2 + \| \lambda x\|^2 \bigg\}
\end{equation}
with $\lambda \ge 0$. The minimizing $x$ has a closed form solution of 
\begin{equation} \label{eq:RR_solution}
    \hat x_{\text{RR}_{\lambda}} = (\A^T \A + \lambda^2 I)^{-1} \A^T b.
\end{equation}
There is a structural connection between the RR and TLS; their solutions are nearly identical with TLS subtracting from the diagonal of $\A^T \A$ and RR adding to it. This can be understood from a Bayesian statistics point of view. Large values of $\lambda$ correspond to stronger evidence of lower variance in $x$ (and zero mean). In the case of an uncertain data matrix, we should have \textit{less} confidence in low variance. 

Golub, Hansen, and O'Leary explored the interplay between TLS and RR in \cite{golub1999tikhonov}. They considered a generalized version of the regularized total least squares (RRTLS) problem 
\begin{align} \label{eq:RRTLS}
\min_{x,\dA,\eta}  \quad &\big\| [\dA,  \ \eta] \big\|_F  \nonumber \\
\text{subject to } &(\A + \dA)x = b + \eta \ \  \text{and} \ \  \|x\| \le \gamma. 
\end{align}
The second constraint is equivalent to imposing a regularization term on the objective of the TLS problem. When the inequality constraint is replaced by an equality, the solution to \eqref{eq:RRTLS}, $\hat x_{\text{RRTLS}_{\alpha}}$ say, is given by the $x$ that solves
\begin{equation}\label{eq:RRTLS_solution}
	x = \left(A^T A + \alpha I  \right)^{-1} A^T b \ \  \text{with} \ \  \alpha = \mu (1+\gamma^2) - \frac{\|b-A  x  \|^2}{1+\gamma^2}
\end{equation}
and $\mu$ the corresponding Lagrange multiplier. Note that for all $\gamma >0$, we have $\alpha \ge -\sigma_{n+1}^2$. Similarly, $\alpha = -\sigma_{n+1}^2$ when $\gamma = \|\hat x_{\text{TLS}}\|$ and $\alpha = 0$ when $\gamma = \|\hat x_{\text{OLS}} \|$. This suggests that regularized total least squares merely adds a standard regularization term to the poorly conditioned TLS matrix, $A^T A - \sigma_{n+1}^2 I$.

Given the interplay between TLS, RR, and RRTLS, it is reasonable to solve either TLS or RR since they are effectively the same problem with different parameters; the additional term becomes regularizing when $\alpha = \lambda^2 > 0$. Although TLS is most appropriate for the case of an uncertain matrix, its ``deregularizing'' effect inflames issues with conditioning making it an unpopular choice for solving typical linear inverse problems. 

There is also the tricky consideration of choosing an appropriate regularization parameter, whether in the standard RR form of \eqref{eq:RR} or $\gamma$ in \eqref{eq:RRTLS}. Although there are many approaches for choosing a parameter such as Morozov's discrepancy principle (MDP), the unbiased predictive risk estimator method (UPR), the generalized cross validation method (GCV), or the ``elbow'' or ``L'' method to name a few, parameter selection is based on heuristic arguments or requires unknown information \textit{a priori}. A detailed treatment of the above methods and their analysis can be found in \cite{vogel2002computational}. Through the remainder of this paper, we drop our discussion of RRTLS, focusing instead on TLS and RR.

\subsection{Robust Least Squares}
The central focus of robust optimization (RO) is to find a solution that is feasible over all possible realizations of uncertain variables. In our case, $\dA \in \U$ is unknown where $\U = \{M \in \R^{m \times n}: |M| \le D \}$ and $D$ is an element-wise positive matrix chosen in a principled fashion, i.e., the degree to which matrix entries are quantized. Our robustified version is written as a minimax optimization problem in (\ref{eq:main}) with the appropriate $\U$. The inner maximization problem guards against over-fitting, effectively regularizing our solution, and reducing sensitivity to perturbations in $\A$.
Note that RO avoids placing a statistical prior on $\Delta$, which can be a strength or weakness depending on the model.

The outline of the rest of the paper is as follows: in section II we derive a closed form solution to the inner maximization problem thereby showing its tractability, section III discusses computational methods for solving problem (\ref{eq:main}), and section IV provides the results to numerical experiments.

\section{Closed Form Solution of Inner Objective and Theory}
\subsection{Floating Point Uncertainty}
We begin by breaking \eqref{eq:main} into an inner maximization and outer minimization problem:
\begin{equation} \label{eq:floatingPtOptimizationProblem}
\min_x \, \underbrace{ \left\{ \max_{ |\dA| \le D }\;  \| (\A+\dA)x -b \|^2 \right\} }_{f(x)} 
\end{equation}
which is equivalent to
\begin{equation}
\min_x   f(x) \quad \text{subject to} \quad  f(x) = \max_{|\dA| \le D} \big\|(A+\dA)x - b \big\|^2.
\end{equation}
Here $|\dA| \le D$ indicates that the magnitude of elements in $\dA$ are bound by the corresponding non-negative components of $D$. That is, $D_{i,j}$ constrains element $\dA_{i,j}$ of the uncertainty matrix. The inner maximization problem will recover $f(x)$. Because the function $x \mapsto \| (\A+\dA)x -b \|^2$ is convex in $x$ and the supremum over an arbitrary family of convex functions is convex, it follows that $f$ is convex making $\min_x f(x)$ an unconstrained convex optimization problem. To evaluate $f$ and find a subgradient, it remains to find the maximizing $\dA \in \R^{m \times n}$. 
\begin{theorem} \label{thm:floating}
	The maximizing function $f(x)$ in \eqref{eq:floatingPtOptimizationProblem} is given by 
	\begin{equation}
		 f(x) = \big\|Ax-b \big\|^2 + 2 \langle \, |Ax-b|, D|x| \, \rangle + \big\| \, D|x| \, \big\|^2
	\end{equation}
	where $|x|$ and $|Ax-b|$ denote the vectors of component-wise absolute values. 
	\begin{proof}
		Since $f(x) = \max_{|\dA| \le D} \big\|(A+\dA)x - b \big\|^2$ must be maximized over $\Delta$, we can fix $x$, define $c:= \A x - b$, then treat it as constant. Exploiting row-wise separability, we write 
		\begin{equation} \label{eq:rowwisemax}
		\max_{|\dA| \le D} \ \|\dA x + c \|^2 = \sum_{i = 1}^m \max_{|\dAi^T| \le D_i^T}  (\dAi^T x + c_i)^2
		\end{equation}  
		where $\dAi^T$  and $D_i^T$ are the $i\ith$ row of $\dA$ and $D$, respectively, and $c_i$ is the $i\ith$ element of vector $c$. We now work row by row. Note that we can switch to absolute values rather than squares when maximizing for each row. Applying the triangle inequality gives an upper bound
		\begin{align} \label{eq:upperbound}
		|\dAi^T x + c_i| \  \le & \   |c_i| +  \sum_{j = 1}^n |\dA_{i,j}| \, | x_j| \\
		 \le & \  |c_i| + \sum_{j = 1}^n D_{i,j} |x_j| \quad = \quad D_i^T 
		|x|  + |c_i|. \nonumber
		\end{align}
		It is easily verified that the upper bound is achieved when
		\begin{equation} \label{eq:optimaldelta}
		\dA = D \odot \sign(x \, c^T)
		\end{equation}
		making it a solution to the inner maximization problem. Here, $\odot$ denotes the Hadamard or element-wise product. Recalling that $c_i = [Ax - b]_i$, we simplify to 
		\begin{equation}
			f(x) =  \big\|Ax-b \big\|^2 + 2 \langle \, |Ax-b|, D|x| \, \rangle + \big\| \, D|x| \, \big\|^2.
		\end{equation}
	\end{proof}
\end{theorem}
Using Theorem \ref{thm:floating}, the optimization problem in \eqref{eq:floatingPtOptimizationProblem} can be rewritten as 		
\begin{equation} \label{eq:closedFormFloatingPoint}
\min_x f(x) = \min_x \bigg\{ \big\|Ax-b \big\|^2 + 2 \langle \, |Ax-b|, D|x| \, \rangle + \big\| \, D|x| \, \big\|^2      \bigg\}.
\end{equation}
In general, $f$ is not differentiable $\big($e.g., let $m=n=1$ and $A=D=b=1$, then $f(x)=(x-1)^2 + 2|x|\cdot|x-1| + x^2$ is not differentiable at $x\in\{0,1\}\big)$, but we are guaranteed a subgradient by virtue of convexity. Furthermore, a generalization of Danskin's Theorem~\cite{bertsThesis} provides us with a method to find elements of the subgradient for all $x$. In particular,
\begin{equation} \label{eq:subgradient}
f'(x) := 2 (\A+\dA_x)^T\big[(\A+\dA_x)x - b \big] \  \in \ \partial f(x)
\end{equation}
where $\dA_x$ indicates the optimal $\dA$ for a given $x$ as provided in Eq.~\eqref{eq:optimaldelta}.

\subsection{Fixed Point Uncertainty}
Fixed point uncertainty is a special case of \eqref{eq:floatingPtOptimizationProblem} with all elements bound by the same value. We can write this constraint as $\| \Delta \|_\infty \le \delta$ with $\| \dA \|_{\infty}$ representing the largest magnitude element of $\dA$. Problem \eqref{eq:floatingPtOptimizationProblem} becomes
\begin{equation} \label{eq:Robust}
\min_x \, \left\{  \max_\constraint \; \| (\A+\dA)x -b \|^2 \right\}
\end{equation}
with a solution denoted by $\hat x_{\text{RO}}$. We focus on fixed point error for the remainder of the paper. This instance gives rise to the following corollary.

\begin{corollary}
	For fixed point uncertainty, problem \eqref{eq:Robust} reduces to 
	\begin{equation}\label{eq:robustrecast}
		\min_x \bigg\{ \|\A x-b\|^2 + 2 \delta \|x\|_1 \|\A x-b\|_1 + m \delta^2 \|x\|_1^2 \bigg\}.
	\end{equation}
	\begin{proof}
		We take $\oner_k \in \R^k$ to be the vector composed of ones. Note that $\| \Delta \|_\infty \le \delta$ is equivalent to $|\Delta| \le D$ when $D =  \delta \, \oner_m^{}  \oner_n^T$. By Theorem \ref{thm:floating}, the result follows easily using properties of inner products.
	\end{proof}
\end{corollary}

\subsection{Explicit Regularization of Robust Objective}
The robust formulation presented above is intended to address uncertainty in the data matrix from rounding error. Although implicit regularization occurs via the 3rd term in both \eqref{eq:closedFormFloatingPoint} and \eqref{eq:robustrecast}, we don't address poor conditioning of $A$ directly. To illustrate this point, note that $\hat x_{\text{RO}} \rightarrow \hat x_{\text{OLS}}$ as $\delta \rightarrow 0$. If $A$ is poorly conditioned, such behavior is undesirable. This can be remedied by adding an $\ell_2$ regularization term and solving
\begin{small}
\begin{equation} \label{eq:regularized_robust_problem}
\min_x \bigg\{   \left( \|\A x-b\|^2 + 2 \delta \|x\|_1 \|\A x-b\|_1 + m \delta^2 \|x\|_1^2 \right) + \lambda^2 \|x\|^2      \bigg\}.
\end{equation}
\end{small}
Since the robust objective and regularization term are both convex, so is their sum. Furthermore, the modified objective is $2 \lambda^2$ strongly convex yielding a unique minimizer, $\hat x_{\text{RRO}_\lambda}$. Strong convexity also implies a bounded sequence of iterates, i.e., $\|x_k - x^*\| \le M$,  and therefore guarantees convergence under the mirror-descent method. In fact, the same technique works with regularizers other than $\ell_2$, and can be solved via proximal mirror-descent methods when the regularizer is smooth or has an easy-to-compute proximity operator~\cite[Ch. 9]{BeckBook2017}.

\section{Algorithms} \label{section:algos}
Equipped with an element of the subdifferential, we can employ a variety of solvers. Bundle methods are a promising choice as they sequentially form an approximation to the objective. This is done by using the location, function value, and subgradient of previous iterates to lower bound the objective with supporting hyperplanes. At each step, a direction finding quadratic program must be solved, but can be dealt with rapidly thanks to software such as CVXGEN \cite{cvxgen} and ECOS \cite{domahidi2013ecos}. A survey on the topic can be found in \cite{makela2002survey}.

Since the problem is unconstrained, another option is to use smooth optimization techniques, and in particular quasi-Newton methods that form low-rank approximations of $\nabla^2 f$. Our objective $f$ is not differentiable, much less twice so, and consequently methods requiring second derivatives seem theoretically inappropriate. However, there is a growing body of literature going back to Lemar\'echal \cite{32} recognizing the empirical success of these methods. See \cite{29,31} and references therein. The success of these methods depends on the smoothness near the solution. Using the \verb|minFunc| MATLAB package \cite{SchmidtCode} with random and zero vector initializations, and using the package's default solver of limited-memory BFGS~\cite{LBFGS}, we observed fast and accurate convergence for several test problems. We used this method for the numerical experiments presented in the next section.

Subgradient descent is an appealing option due to its simplicity and flexibility. The method can easily handle constraints via projections or a regularization term $h(x)$ with proximal operators. Recent advances in proximal subgradient descent methods and their analysis can be found in \cite{cruz2017proximal} and \cite{millan2018inexact}. Convergence may be slow, but is less of a concern with heavy quantization. The 3rd terms in \eqref{eq:closedFormFloatingPoint} and \eqref{eq:robustrecast} effectively regularize the objective when $\delta$ and $\|D\|$ are large yielding faster convergence. A drawback of the subgradient descent method is that a step-size scheme must be chosen in advance. One popular choice is the diminishing, non-summable step length given by $t_k = \frac{1}{\sqrt{k+1} \cdot \|f'(x_k)\|}$. A discussion on the choice of step lengths can be found in \cite{boyd2003subgradient}. Stopping criteria are also difficult to specify since it is not a true descent method (the objective is not guaranteed to decrease on every iteration). The algorithm therefore typically runs for a fixed number of iterations.

\section{Numerical Experiments} \label{section:experiments}
\paragraph*{\textbf{Setup}}We use MATLAB's pseudo-random number generator \verb|randn| to draw 10,000 standard normal matrices, $\bar A~\in~\R^{30 \times 15}$. The bar notation in this section indicates true and unobserved values. We fixed the condition number of our matrices at 100. We did so by performing a singular value decomposition (SVD) such that $\bar A = U \Sigma V^T$ where $U$, $V$ are unitary and $\Sigma$ is diagonal, then replaced the diagonal of $\Sigma$ by linearly decaying values from $1 \rightarrow \frac{1}{100}$. This ensures our test matrices show sensitivity to noisy measurements prior to quantization.

\begin{figure}  
	\begin{minipage}{.5\textwidth}
		\centering
		\captionsetup{width=.9\linewidth}
		\includegraphics[width=.9\linewidth]{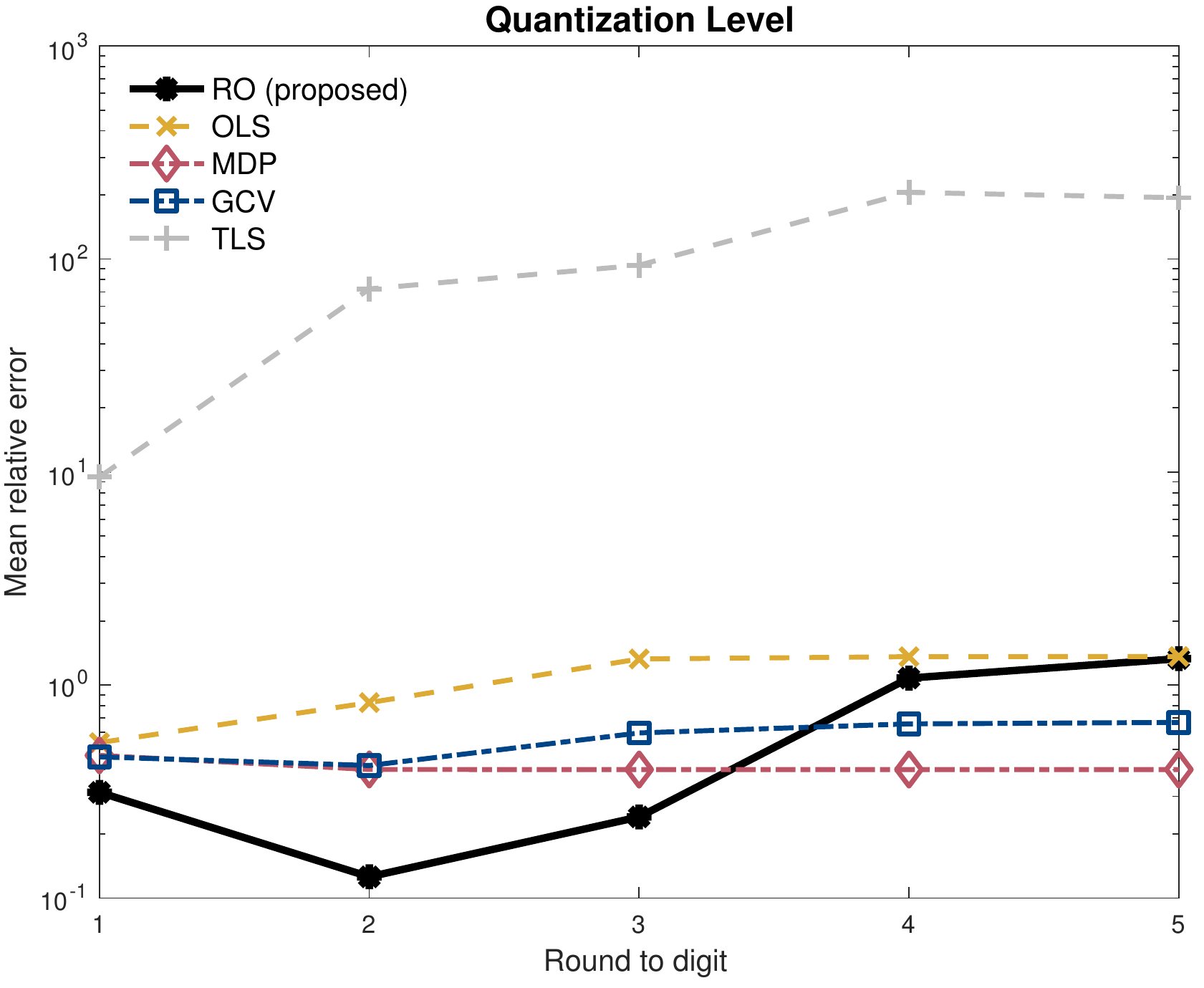}
		\caption{Mean relative error $\|e\| = \|\hat x - \bar x\|/\|\bar x\|$ over 10k simulations as a function of the digit to which matrix entries were rounded to. Note that $\delta = 0.5\times 10^{-\text{Round to digit}}.$  }
		\label{figs:errors}
	\end{minipage}
	\begin{minipage}{.5\textwidth}
		\centering
		\captionsetup{width=.9\linewidth}
		\includegraphics[width=.9\linewidth]{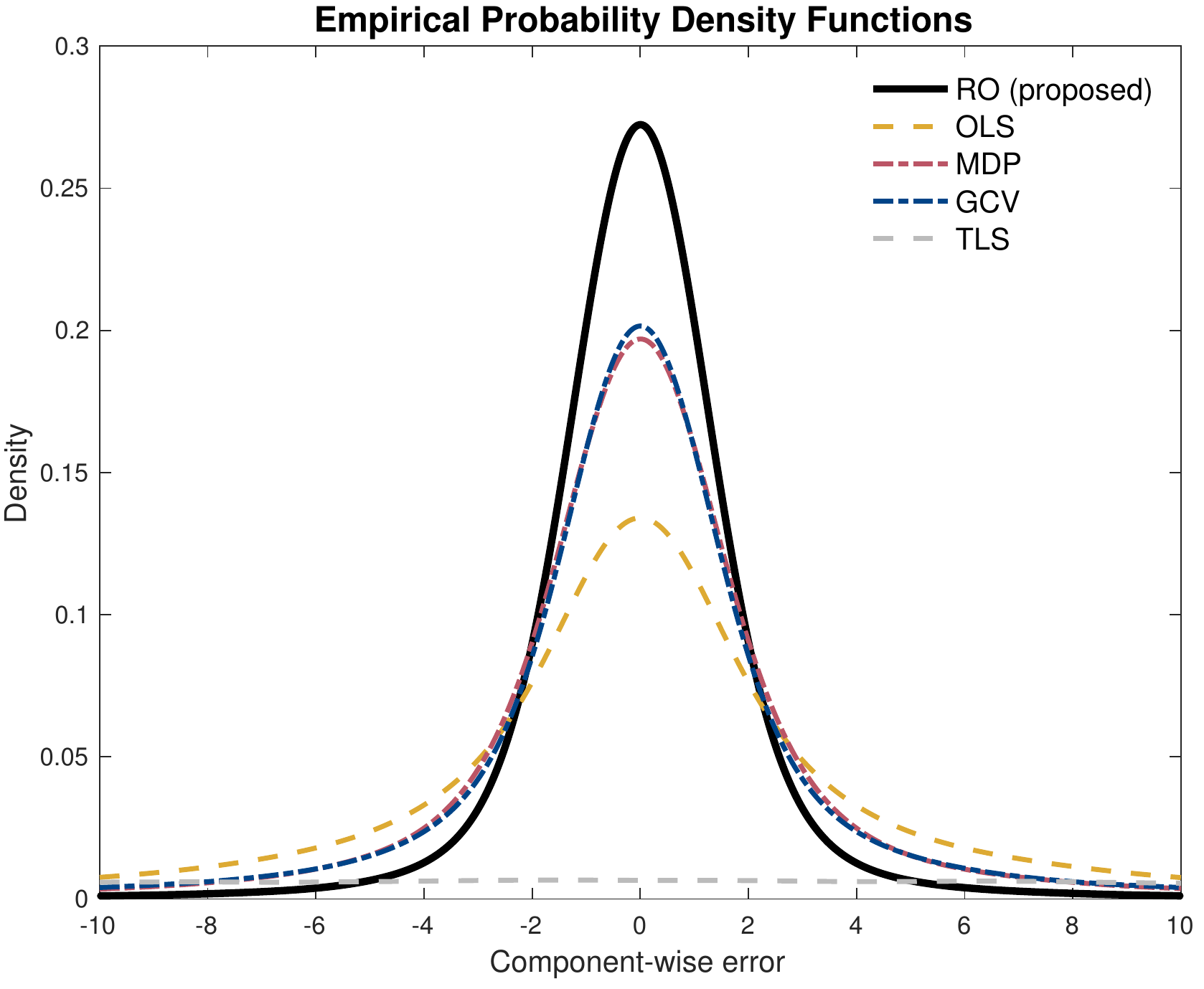}
		\caption{Empirical probability density function of component-wise error over 10k simulation when $\bar x$ is drawn from a Cauchy distribution and $\bar A$ is quantized to the hundredth spot (round to digit = 2 and   $ \delta = \num{0.5e-2}$ ).}
		\label{figs:cauchy_density}
	\end{minipage}
\end{figure}

We conducted two experiments. In the first, we draw random vectors from a heavy-tailed Cauchy distribution with median $0$ and scale parameter $1$. In the second, the true solution $\bar x$ has its first element drawn from $\bar x_1 = \pm 100$ with a sign drawn uniformly at random and the remaining elements drawn from a standard normal distribution, $\bar x_{2:15} \sim \mathcal{N}(0,I)$. This setup exposes the bias of RR solutions. We remark that we do this for illustrative purposes; if one suspected a signal of having such a large element, it would of course make sense to run an outlier detection method first.

We obtain our ``true'' right hand side by taking the image of $\bar x$ under $\bar A$, that is, $\bar b := \bar A \bar x$, then generate our observed measurement vector by letting $b := \bar b + \eta$ where $\eta \sim \mathcal{N}\left(0, \frac{\|\bar b\|^2}{m \cdot \text{SNR} } I\right)$ with the desired signal-to-noise ratio (SNR) fixed at $\text{SNR} = 50$. Finally, we quantize $\bar A$ which yields our observed $A$, then solve problems for OLS~\eqref{eq:OLS}, TLS~\eqref{eq:TLS}, $\text{RR}_\lambda$~\eqref{eq:RR}, RO~\eqref{eq:Robust}, and $\text{RRO}_\lambda$~\eqref{eq:regularized_robust_problem}. The $\lambda$ for $\text{RR}_\lambda$ and $\text{RRO}_\lambda$ are chosen according to the GCV and MDP criteria discussed in Section~\ref{sec:Tikhonov}, and the $\delta$ for $\text{RO}$ and $\text{RRO}_\lambda$ is based on the observed accuracy, e.g., if $\bar A$ is rounded to the $2^{\text{nd}}$ decimal, then  $\delta = 0.5 \cdot 10^{-2}$. We use the \verb|minfunc| implementation of limited-memory MFGS to solve the robust problems with a random initialization. The GCV parameters is recovered using an approximate line search for the corresponding objectives provided in Vogel's text \cite{vogel2002computational}. We omit UPR because of its similarity to the GCV parameter.

The MDP parameter is found using the bisection method to solve $\|A \, \hat x_{\text{RR}_\lambda} - b\|^2 - \rho = 0$ for $\lambda$, where $\hat x_{\text{RR}_\lambda}$ is the solution to \eqref{eq:RR} and the parameter $\rho$ ideally reflects the residual's noise floor. We chose this method since it is computationally cheap and easy to program, though it is also possible to solve a constrained least-squares problem directly using standard software like \verb|cvxopt| \cite{cvxopt}. 
The goal is to find $\lambda$ such that $\|A \hat x_{\text{RR}_\lambda}- b\|^2$ equals the expected uncertainty introduced by noise.
For $\bar A \bar x - b = \eta \sim \mathcal{N}(0, \sigma^2 I)$ which is an $m$ component Gaussian random variable, 
\begin{equation}
    \frac{1}{\sigma^2} \|\eta\|^2 = \sum_{i=1}^m \left( \frac{\eta_i}{\sigma} \right)^2 \sim \chi^2(m)
\end{equation} 
since $\frac{\eta_i}{\sigma} \sim \mathcal{N}(0,1)$ and $\left( \frac{\eta_i}{\sigma}\right)^2 \sim \chi^2(1)$. When $m=30$, as in our case, both the mean and median are approximately 30 (exact for mean). To ensure feasibility, i.e., a real root, we choose $\rho$ such that $\mathbb{P}\left(\|\bar A \bar x-b\|^2 - \rho \le 0 \right) = 95\%$, which can be easily calculated using $\chi^2$ inverse CDF tables. For our experiments, $\rho \approx \frac{2 \|\bar b \|^2}{3 \cdot \text{SNR}}$. An appropriate $\rho$ is generally unknown \textit{a priori} unless one knows the noise variance.

\begin{figure}  
	\begin{minipage}{0.5\textwidth}
		\centering
		\captionsetup{width=.9\linewidth}
		\includegraphics[width=.9 \textwidth]{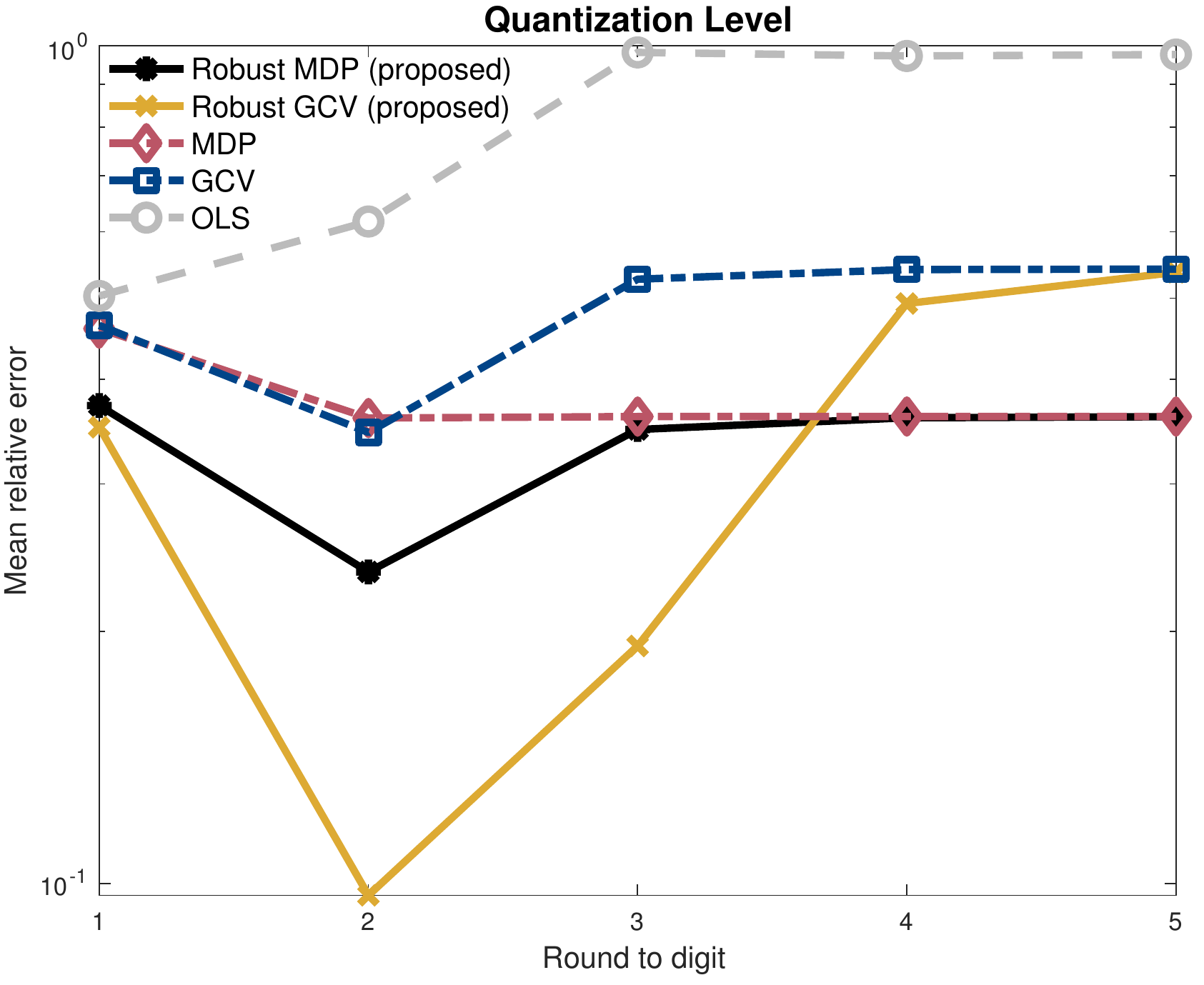} 
		\caption{Mean relative error $\|e\| = \|\hat x - \bar x\| / \|\bar x\|$ over 10k simulations of $\text{RR}_\lambda$ and its robust counter part,  $\text{RRO}_\lambda$. Note that $\delta = 0.5\times 10^{-\text{Round to digit}}$. The robust counterpart shows benefit for large quantization then converges to corresponding RR solution.} 
		\label{figs:regularized_mean_error}
	\end{minipage}
	\begin{minipage}{0.5\textwidth}
		\centering
		\captionsetup{width=.9\linewidth}
		\includegraphics[width=.9 \textwidth]{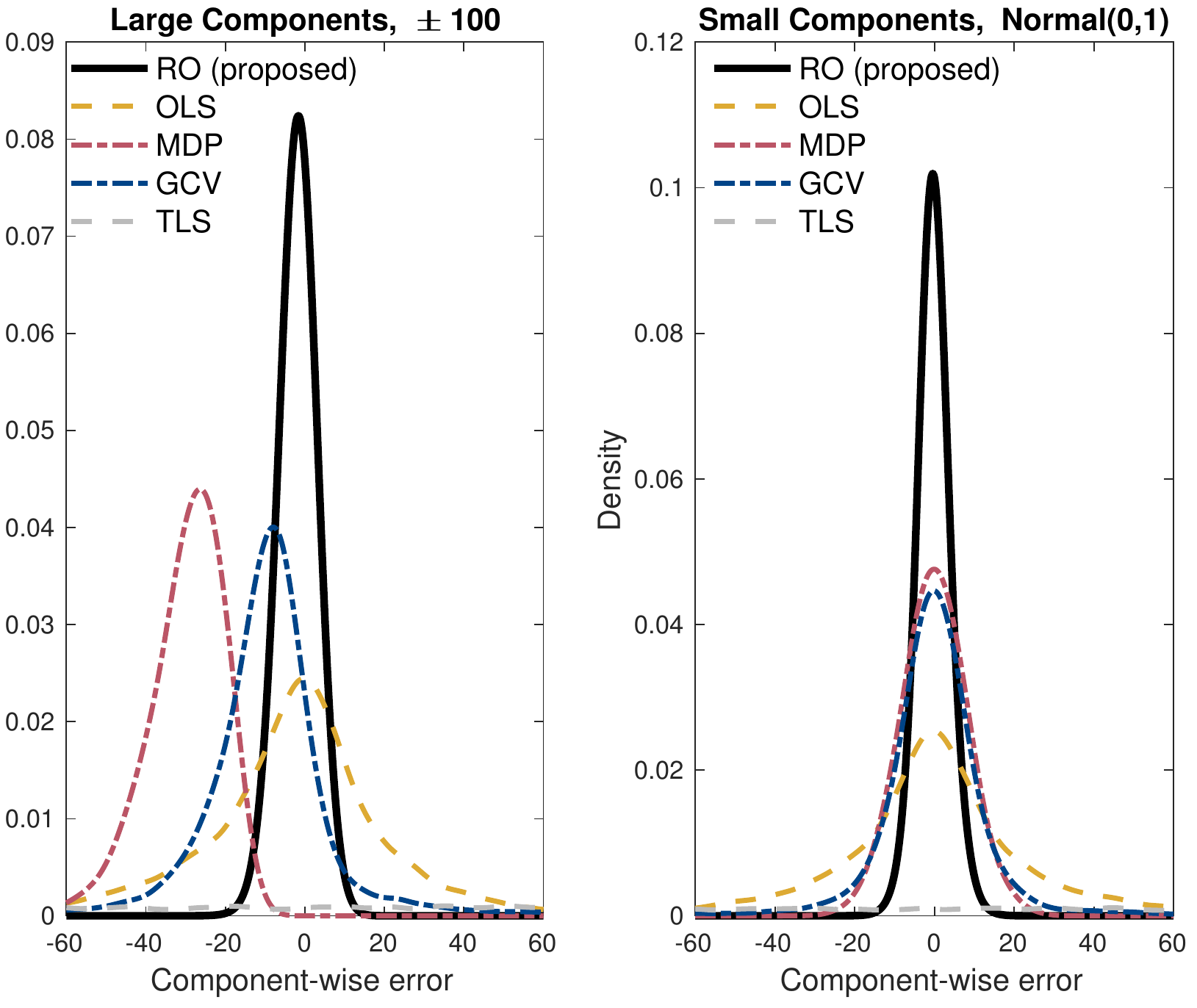}
		\caption{Empirical probability densities of  component-wise errors for the single large element experiment when $\bar A$ is quantized to the hundredth. Left: Pdf of the error (adjusted by the sign of the large element $\bar{x}_1$)
			of different methods for the large element ($\bar{x}_1 = \pm 100$). Note that RR and GCV have have modes away from zero indicating bias. Right: Pdf of the error for the  $2^\text{nd}$--$15^\text{th}$ components (which are drawn from $\mathcal{N}(0,1)$).} 
		\label{figs:onelarge_densities}
	\end{minipage}
\end{figure}

\paragraph*{\textbf{Numerical Results}}
Results for the Cauchy distributed measurement vector are displayed in Figure~\ref{figs:errors}, which shows mean error as a function of rounding digit for $\bar A$, and Figure~\ref{figs:cauchy_density}, which illustrates an empirical probability density function (pdf) of component-wise errors for different methods using a standard Gaussian kernel. In this case, $A$ is quantized to the hundredth spot. Note that the robust solution does well at reducing error in the face of heavy quantization. This is evidenced by low mean error through the thousandth spot and small variance as seen in the pdf.
The RO solution converges to the OLS solution as expected when $\delta \rightarrow 0$, i.e., when the quantization effect is small. This might be undesirable, especially when $A$ is ill-conditioned, and can be addressed via the inclusion of a regularization term. The mean relative error performance for the regularized robust problem given by \eqref{eq:regularized_robust_problem} is shown in Figure \ref{figs:regularized_mean_error}. Since the proposed method does not improve estimation universally, we recommend its use when quantization results in precision loss of greater than $0.1\%$. We reiterate that when quantization effects are heavy, there is a notable benefit of the proposed method and the parameter $\delta$ can be chosen in a natural way. 

A second observation is that RO doesn't sacrifice accuracy for bias as RR does. This phenomenon can be observed in Figure \ref{figs:onelarge_densities} which depicts empirical pdf's of component- wise error for different methods. The plots to the right show error behavior for small components of the solution, that is, for $\bar x_{2:15} \sim \mathcal{N}(0,I)$. Since values are close to zero, RR is expected to perform well at estimating the true solution. Indeed, RR implicitly assumes a prior with mean zero. 
The densities on the left show errors for ``large'' components drawn from $\{\pm 100\}$ (and the sign of the error is adjusted by the sign of the large component, so a negative error indicates that the estimate is biased toward zero, and a positive error indicates bias to $\pm\infty$). Rather than localizing about zero, the absolute error's mode is observed around -7 for GCV and UPR and -25 for MDP. Heavier penalty terms (bigger $\lambda$) place less weight on the LS term in the objective and bias estimates more aggressively towards zero. Typical $\lambda$ values for GCV, UPR, and MDP are 0.041, 0.042, and 0.17, respectively. The robust solution accurately estimates the large and small components values without the extreme errors observed with OLS and TLS; that is, RO appears to have a low-variance, in-line with RR, but also low-bias, comparable or better than OLS.

\section{Conclusion}
In this paper, we presented a robust method for addressing fixed and floating point uncertainty due to quantization error. After reviewing the limitations of existing methods, we formulated a tractable robust objective and presented algorithms to solve it. The only parameter necessary in our formulation is chosen in a principled fashion, i.e., by observing the degree to which matrix elements are rounded. Our numerical experiments show that robust least squares outperforms OLS, TLS, and RR, effectively balancing bias and accuracy. 

\subsection*{Acknowledgments}
We would like to thank Marek Petrik for his keen eye and helping us spot an elegant simplification to our approach.  

\vspace{.25in}

\bibliographystyle{ieeetr}

\bibliography{thesisBecker}

\end{document}